\documentclass[journal]{IEEEtran}

\usepackage{subfigure}
\usepackage{bm}
\usepackage{amssymb}
\usepackage{epsfig}
\usepackage{amsmath}
\usepackage{amsthm}
\usepackage{xcolor}
\usepackage{graphicx}
\usepackage{epstopdf}
\usepackage{amsfonts}%
\usepackage{indentfirst}
\usepackage{color}
\newtheorem{problem}{\textbf{Problem}}
\newtheorem{definition}{\textbf{Definition}}

\newtheorem{theorem}{\rm\textbf{Theorem}}
\newtheorem{lemma}{\rm\textbf{Lemma}}

\newtheorem{remark}{\rm\textbf{Remark}}

\begin{document}
%
\title{Sufficient Conditions for Feasibility of Optimal Control Problems Using Control Barrier Functions}
%
%
%

\author{Wei Xiao,~\IEEEmembership{}
        Calin Belta,~\IEEEmembership{}
          and Christos G. Cassandras~\IEEEmembership{}
\thanks{The authors are
	with the Division of Systems Engineering and Center for Information and
	Systems Engineering, Boston University, Brookline, MA, 02446, USA
	\texttt{{\small \{xiaowei,cbelta,cgc\}@bu.edu}}}
\thanks{This work was
	supported in part by NSF under grants IIS-1723995, CPS-1446151, ECCS-1931600,
	DMS-1664644, CNS-1645681, by AFOSR under grant FA9550-19-1-0158, by ARPA-E's
	NEXTCAR program under grant DE-AR0000796 and by the MathWorks. } 
}

\maketitle

\begin{abstract}
It has been shown that satisfying state and control
constraints  while optimizing quadratic costs subject to desired (sets of)
state convergence for affine control systems can be reduced to a sequence of quadratic programs (QPs) by using
Control Barrier Functions (CBFs) and Control Lyapunov Functions (CLFs). One of the main challenges in this approach is ensuring the feasibility of these QPs, especially under tight control bounds and safety constraints of high relative degree. In this paper, we provide sufficient conditions for guranteed feasibility. The sufficient conditions are captured by a single constraint that is enforced by a CBF, which is added to the QPs such that their feasibility is always guaranteed. The additional constraint is designed to be always compatible with the existing constraints, therefore, it cannot make a feasible set of constraints infeasible - it can only increase the overall feasibility. We illustrate the effectiveness of the proposed approach on an adaptive cruise control problem.
\end{abstract}


%
\IEEEpeerreviewmaketitle

\section{INTRODUCTION}

\label{sec:intro}

Constrained optimal control problems with safety specifications are central to increasingly widespread safety critical autonomous and cyber physical systems. Traditional Hamiltonian analysis \cite{Bryson1969} and dynamic programming \cite{Denardo2003} cannot accommodate the size and nonlinearities of such systems, and their applicability is mostly limited to linear systems. 
Model Predictive Control (MPC) \cite{Rawlings2018} methods have been shown to work for large, non-linear systems. However, safety requirements are hard to be guaranteed between time intervals in MPC. Motivated by these limitations, barrier and control barrier functions enforcing safety have received increased attention in the past years \cite{Aaron2014} \cite{Glotfelter2017} \cite{Xiao2019}. 

Barrier functions (BFs) are Lyapunov-like functions \cite{Tee2009},
\cite{Wieland2007}, whose use can be traced back to optimization problems
\cite{Boyd2004}. More recently, they have been employed to prove set
invariance \cite{Aubin2009}, \cite{Prajna2007}, \cite{Wisniewski2013} and to address
multi-objective control problems \cite{Panagou2013}. In \cite{Tee2009}, it was proved
that if a BF for a given set satisfies Lyapunov-like conditions, then the set
is forward invariant. A less restrictive form of a BF, which is allowed to
grow when far away from the boundary of the set, was proposed in
\cite{Aaron2014}. Another approach that allows a BF to be zero was proposed in
\cite{Glotfelter2017}, \cite{Lindemann2018}. This simpler form has also been
considered in time-varying cases and applied to enforce Signal Temporal Logic
(STL) formulas as hard constraints \cite{Lindemann2018}.

Control BFs (CBFs) are extensions of BFs for control systems, and are used to
map a constraint defined over system states to a constraint on the control
input. The CBFs from \cite{Aaron2014} and \cite{Glotfelter2017} work for constraints
that have relative degree one with respect to the system dynamics. A
backstepping approach was introduced in \cite{Hsu2015} to address higher
relative degree constraints, and it was shown to work for relative degree two.
A CBF method for position-based constraints with relative degree two was also
proposed in \cite{Wu2015}. A more general form was considered in \cite{Nguyen2016}, which works
for arbitrarily high relative degree constraints, employs input-output
linearization and finds a pole placement controller with negative poles to
stabilize an exponential CBF to zero. The high order
CBF (HOCBF) proposed in \cite{Xiao2019} is simpler and more general than the
exponential CBF \cite{Nguyen2016}. 

Most works using CBFs to enforce safety are based on the assumption that the (nonlinear) control system is affine in controls and the cost is quadratic in controls. Convergence to desired states is achieved by using Control Lyapunov Functions (CLFs) \cite{Aaron2012}.  The time domain is discretized, and the state is assumed to be constant within each time step (at its value at the beginning of the step). The optimal control problem becomes a Quadratic Program (QP) in each time step, and the optimal control value is kept constant over each such step. Using this approach, the original optimal control problem is reduced to a (possibly large) sequence of quadratic programs (QP) - one for each interval \cite{Galloway2013}. While computationally efficient, this myopic approach can easily lead to infeasibility: the constant optimal control derived at the beginning of an interval can lead the system to a state that gives incompatible control constraints  at the end of the interval, rendering the QP corresponding to the next time interval infeasible.     

For the particular case of an adaptive cruise control (ACC) problem in \cite{Aaron2014}, it was shown that an additional constraint (minimum braking distance) can help keep the system away from states 
leading to incompatibility of control CBF and CLF constraints. However, this additional constraint itself may conflict with other constraints in the ACC problem, such as the control bounds. To improve the problem feasibility for general optimal control problems with the CBF method, the penalty method \cite{Xiao2019} and adaptive CBF \cite{Xiao2020} were proposed; however, they still do not guarantee the QP feasibility. 

In this paper, we provide a novel method to find sufficient conditions to guarantee the feasibility of CBF-CLF based QPs. This is achieved by the proposed feasibility constraint method that makes the problem constraints compatible in terms of control given an arbitrary system state. The sufficient conditions are captured by a single constraint that is enforced by a CBF, and is added to the problem to formulate the sequence of QPs mentioned above with guaranteed feasibility. The added constraint is always compatible with the existing constraints and, therefore, it cannot make a feasible set of constraints infeasible. However, by ``shaping'' the constraint set of a current QP, it guarantees the feasibility of the next QP in the sequence. We illustrate our approach and compare it to other methods on an ACC problem.

The remainder of the paper is organized as follows. In Sec.\ref{sec:pre}, we provide preliminaries on HOCBF and CLF. Sec.\ref{sec:prob} formulates an optimal control problem and outlines our CBF-based solution approach. We show how we can find a feasibility constraint for an optimal control problem in Sec.\ref{sec:fc}, and present case studies and simulation results in Sec. \ref{sec:case}. We conclude the paper in Sec.\ref{sec:conclusion}.

\section{PRELIMINARIES}

\label{sec:pre}

\begin{definition}
	\label{def:classk} (\textit{Class $\mathcal{K}$ function} \cite{Khalil2002}) A
	continuous function $\alpha:[0,a)\rightarrow[0,\infty), a > 0$ is said to
	belong to class $\mathcal{K}$ if it is strictly increasing and $\alpha(0)=0$.
\end{definition}

Consider an affine control system of the form
\begin{equation}
\dot{\bm{x}}=f(\bm x)+g(\bm x)\bm u \label{eqn:affine}%
\end{equation}
where $\bm x\in X\subset\mathbb{R}^{n}$, $f:\mathbb{R}^{n}\rightarrow\mathbb{R}^{n}$
and $g:\mathbb{R}^{n}\rightarrow\mathbb{R}^{n\times q}$ are locally
Lipschitz, and $\bm u\in U\subset\mathbb{R}^{q}$ is the control constraint set
defined as
\begin{equation}
U:=\{\bm u\in\mathbb{R}^{q}:\bm u_{min}\leq\bm u\leq\bm u_{max}\}.
\label{eqn:control}%
\end{equation}
with $\bm u_{min},\bm u_{max}\in\mathbb{R}^{q}$ and the inequalities are
interpreted componentwise.

\begin{definition}
	\label{def:forwardinv} A set $C\subset\mathbb{R}^{n}$ is forward invariant for
	system (\ref{eqn:affine}) if its solutions starting at any $\bm x(0) \in C$
	satisfy $\bm x(t)\in C,$ $\forall t\geq0$.
\end{definition}

\begin{definition}
	\label{def:relative} (\textit{Relative degree}) The relative degree of a
	(sufficiently many times) differentiable function $b:\mathbb{R}^{n}%
	\rightarrow\mathbb{R}$ with respect to system (\ref{eqn:affine}) is the number
	of times it needs to be differentiated along its dynamics until the control
	$\bm u$ explicitly shows in the corresponding derivative.
\end{definition}

In this paper, since function $b$ is used to define a constraint $b(\bm
x)\geq0$, we will also refer to the relative degree of $b$ as the relative
degree of the constraint. 

For a constraint $b(\bm x)\geq0$ with relative
degree $m$, $b:\mathbb{R}^{n}\rightarrow\mathbb{R}$, and $\psi_{0}(\bm
x):=b(\bm x)$, we define a sequence of functions $\psi_{i}:\mathbb{R}%
^{n}\rightarrow\mathbb{R},i\in\{1,\dots,m\}$:
\begin{equation}
\begin{aligned} \psi_i(\bm x) := \dot \psi_{i-1}(\bm x) + \alpha_i(\psi_{i-1}(\bm x)),i\in\{1,\dots,m\}, \end{aligned} \label{eqn:functions}%
\end{equation}
where $\alpha_{i}(\cdot),i\in\{1,\dots,m\}$ denotes a $(m-i)^{th}$ order
differentiable class $\mathcal{K}$ function.

We further define a sequence of sets $C_{i}, i\in\{1,\dots,m\}$ associated
with (\ref{eqn:functions}) in the form:
\begin{equation}
\label{eqn:sets}\begin{aligned} C_i := \{\bm x \in \mathbb{R}^n: \psi_{i-1}(\bm x) \geq 0\}, i\in\{1,\dots,m\}. \end{aligned}
\end{equation}

\begin{definition}
	\label{def:hocbf} (\textit{High Order Control Barrier Function (HOCBF)}
	\cite{Xiao2019}) Let $C_{1}, \dots, C_{m}$ be defined by (\ref{eqn:sets}%
	) and $\psi_{1}(\bm x), \dots, \psi_{m}(\bm x)$ be defined by
	(\ref{eqn:functions}). A function $b: \mathbb{R}^{n}\rightarrow\mathbb{R}$ is
	a High Order Control Barrier Function (HOCBF) of relative degree $m$ for
	system (\ref{eqn:affine}) if there exist $(m-i)^{th}$ order differentiable
	class $\mathcal{K}$ functions $\alpha_{i},i\in\{1,\dots,m-1\}$ and a class
	$\mathcal{K}$ function $\alpha_{m}$ such that 
	\begin{equation}
	\label{eqn:constraint}\begin{aligned} \sup_{\bm u\in U}[L_f^{m}b(\bm x) + L_gL_f^{m-1}b(\bm x)\bm u + S(b(\bm x)) \\+ \alpha_m(\psi_{m-1}(\bm x))] \geq 0, \end{aligned}
	\end{equation}
	for all $\bm x\in C_{1} \cap,\dots, \cap C_{m}$. In
	(\ref{eqn:constraint}), $L_{f}^{m}$ ($L_{g}$) denotes Lie derivatives along
	$f$ ($g$) $m$ (one) times, and $S(\cdot)$ denotes the remaining Lie derivatives
	along $f$ with degree less than or equal to $m-1$ (omitted for simplicity, see
	\cite{Xiao2019}).
\end{definition}

The HOCBF is a general form of the relative degree one CBF \cite{Aaron2014},
\cite{Glotfelter2017}, \cite{Lindemann2018} (setting $m=1$ reduces the HOCBF to
the common CBF form in \cite{Aaron2014}, \cite{Glotfelter2017}, \cite{Lindemann2018}), and it is also a general form of the exponential CBF
\cite{Nguyen2016}.

\begin{theorem}
	\label{thm:hocbf} (\cite{Xiao2019}) Given a HOCBF $b(\bm x)$ from Def.
	\ref{def:hocbf} with the associated sets $C_{1}, \dots, C_{m}$ defined
	by (\ref{eqn:sets}), if $\bm x(0) \in C_{1} \cap,\dots,\cap C_{m}$,
	then any Lipschitz continuous controller $\bm u(t)$ that satisfies
	(\ref{eqn:constraint}), $\forall t\geq0$ renders $C_{1}\cap,\dots,
	\cap C_{m}$ forward invariant for system (\ref{eqn:affine}).
\end{theorem}

\begin{definition}
	\label{def:clf} (\textit{Control Lyapunov function (CLF)} \cite{Aaron2012}) A
	continuously differentiable function $V: \mathbb{R}^{n}\rightarrow\mathbb{R}$
	is an exponentially stabilizing control Lyapunov function (CLF) for system
	(\ref{eqn:affine}) if there exist constants $c_{1} >0, c_{2}>0, c_{3}>0$ such
	that for all $\bm x\in X$, $c_{1}||\bm x||^{2} \leq V(\bm x)
	\leq c_{2} ||\bm x||^{2}, $
	\begin{equation} \label{eqn:clf}
	\inf_{\bm u\in U} \lbrack L_{f}V(\bm x)+L_{g}V(\bm x)
	\bm u + c_{3}V(\bm x)\rbrack\leq0.
	\end{equation}
	
\end{definition}

Many existing works \cite{Aaron2014}, \cite{Nguyen2016}, \cite{Yang2019}
combine CBFs for systems with relative degree one with quadratic costs to form
optimization problems. Time is discretized and an optimization problem with
constraints given by the CBFs (inequalities of the form (\ref{eqn:constraint}%
)) is solved at each time step. The inter-sampling effect is considered in \cite{Yang2019}. If convergence to a state is desired, then a
CLF constraint of the form (\ref{eqn:clf}) is added, as in \cite{Aaron2014} \cite{Yang2019}. Note that these
constraints are linear in control since the state value is fixed at the
beginning of the interval, therefore, each optimization problem is a quadratic
program (QP). The optimal control obtained by solving each QP is applied at
the current time step and held constant for the whole interval. The state is
updated using dynamics (\ref{eqn:affine}), and the procedure is repeated.
Replacing CBFs by HOCBFs allows us to handle constraints with arbitrary
relative degree \cite{Xiao2019}. This method works conditioned on the
fact that the QP at every time step is feasible. However, this is not
guaranteed, in particular under tight control bounds. In this paper, we show
how we can find sufficient conditions for the feasibility of the QPs.

\section{PROBLEM FORMULATION AND APPROACH}
\label{sec:prob}

$\textbf{Objective}$: (Minimizing cost) Consider an optimal control problem for the system in (\ref{eqn:affine}) with the cost defined as:
\begin{equation}\label{eqn:cost}
J(\bm u(t)) = \int_{0}^{T}\mathcal{C}(||\bm u(t)||)dt
\end{equation}
where $||\cdot||$ denotes the 2-norm of a vector, $\mathcal{C}(\cdot)$ is a
strictly increasing function of its argument, and $T > 0$. Associated with this problem are the requirements that follow.

$\textbf{State convergence}$: We want the state of system (\ref{eqn:affine}) to reach a point $\bm K\in\mathbb{R}^n$, i.e.,
\begin{equation} \label{eqn:target}
\min_{\bm u(t)}||\bm x(T) - \bm K||^2.
\end{equation}

$\textbf{Constraint 1}$ (Safety constraints): 
System (\ref{eqn:affine}) should always satisfy one or more safety requirements of the form:
\begin{equation} \label{eqn:safetycons}
b(\bm x(t))\geq 0, \forall t\in[0,T].
\end{equation}
where $b: \mathbb{R}^n\rightarrow\mathbb{R}$ is continuously differentiable.

$\textbf{Constraint 2}$ (Control constraints): The control must satisfy (\ref{eqn:control}) for all $t\in[0,T]$.

A control policy for system (\ref{eqn:affine}) is $\bm {feasible}$ if constraints (\ref{eqn:safetycons}) and (\ref{eqn:control}) are satisfied for all times. In this paper, we consider the following problem:

\begin{problem}\label{prob:general}
	Find a feasible control policy for system (\ref{eqn:affine}) such that the cost (\ref{eqn:cost}) is minimized,
	and the state convergence (\ref{eqn:target}) is achieved with the minimum $||\bm x(T) - \bm K||^2$. 
\end{problem}

\textbf{Approach:} We use a HOCBF to enforce (\ref{eqn:safetycons}), and use a relaxed CLF to achieve the convergence requirement (\ref{eqn:target}). If the cost (\ref{eqn:cost}) is quadratic in $\bm u$, then we can formalize Problem \ref{prob:general} using a CBF-CLF-QP approach \cite{Aaron2014}, 
with the CBF replaced by the HOCBF \cite{Xiao2019}: 

\begin{equation}\label{eqn:prob_qp}
\min_{\bm u(t), \delta(t)} \int_{0}^T ||\bm u(t)||^2 + p\delta^2(t) dt
\end{equation}
subject to
\begin{equation}\label{eqn:prob_cbf}
L_f^{m}b(\bm x) + L_gL_f^{m-1}b(\bm x)\bm u + S(b(\bm x)) + \alpha_m(\psi_{m-1}(\bm x)) \geq 0,
\end{equation}
\begin{equation} \label{eqn:prob_clf}
L_{f}V(\bm x)+L_{g}V(\bm x)\bm u + \epsilon V(\bm x)\leq \delta(t),
\end{equation}
\begin{equation} \label{eqn:prob_ctrl}
\bm u_{min}\leq\bm u\leq\bm u_{max},
\end{equation}
where $V(\bm x) = ||\bm x(t) - \bm K||^2$, $c_3 = \epsilon > 0$ in Def. \ref{def:clf}, $p > 0$, and $\delta(t)$ is a relaxation for the CLF constraint. We assume that $b(\bm x)$ has relative degree $m$. The above optimization problem is {\em feasible at a given state} $\bm x$ if all the constraints define a non-empty set for the decision variables $\bm u, \delta$. 

The optimal control problem (\ref{eqn:prob_qp}), (\ref{eqn:prob_cbf}), (\ref{eqn:prob_clf}), (\ref{eqn:prob_ctrl}) with decision variables $\bm u(t), \delta(t)$ 
is usually solved point-wise, as outlined in the end of Sec.\ref{sec:pre}. The time interval $[0,T]$ is divided into a finite number of intervals. At every discrete time $\bar{t}\in[0,T)$ defining the bounds of the intervals, we fix the state $\bm x(\bar{t})$, so that the optimal control problem above becomes a QP.  We obtain an optimal control $\bm u^*(\bar{t})$ and we apply it to system (\ref{eqn:affine}) for the whole interval for which $\bar{t}$ is the lower bound. 

This paper is motivated by the fact that this myopic approach can easily lead to infeasible QPs, especially under tight control bounds. In other words, after we apply the constant $\bm u^*(\bar{t})$ to system (\ref{eqn:affine}) starting at $\bm x(\bar{t})$ for the whole interval that starts at $\bar{t}$, we may end up at a state where the HOCBF constraint (\ref{eqn:prob_cbf}) conflicts with the control bounds (\ref{eqn:prob_ctrl}), which would render the QP corresponding to the next time interval infeasible~\footnote{Note that, since the CLF constraint (\ref{eqn:prob_clf}) is relaxed, it does not affect the feasibility of the QP.}. To avoid this, we define an additional {\em feasibility constraint}:

\begin{definition}\label{def:fc} [{\bf feasibility constraint}]
	Suppose the QP (\ref{eqn:prob_qp}), subject to (\ref{eqn:prob_cbf}), (\ref{eqn:prob_clf}) and (\ref{eqn:prob_ctrl}), is feasible at the current state $\bm x(\bar t), \bar t\in [0,T)$. A constraint $b_F(\bm x)\geq 0$, where $b_F:\mathbb{R}^n\rightarrow \mathbb{R}$, is a feasibility constraint if it makes the QP  corresponding to the next time interval feasible.
\end{definition}

In order to ensure that the QP (\ref{eqn:prob_qp}), subject to (\ref{eqn:prob_cbf}), (\ref{eqn:prob_clf}) and (\ref{eqn:prob_ctrl}), is feasible for the next time interval, a feasibility constraint $b_F(\bm x)\geq 0$ should have two important features: $(i)$ it guarantees that (\ref{eqn:prob_cbf}) and (\ref{eqn:prob_ctrl}) do not conflict, $(ii)$ the feasibility constraint itself does not conflict with both (\ref{eqn:prob_cbf}) and (\ref{eqn:prob_ctrl}) at the same time.

\begin{figure}[thpb]
	\centering
	\includegraphics[scale=0.19]{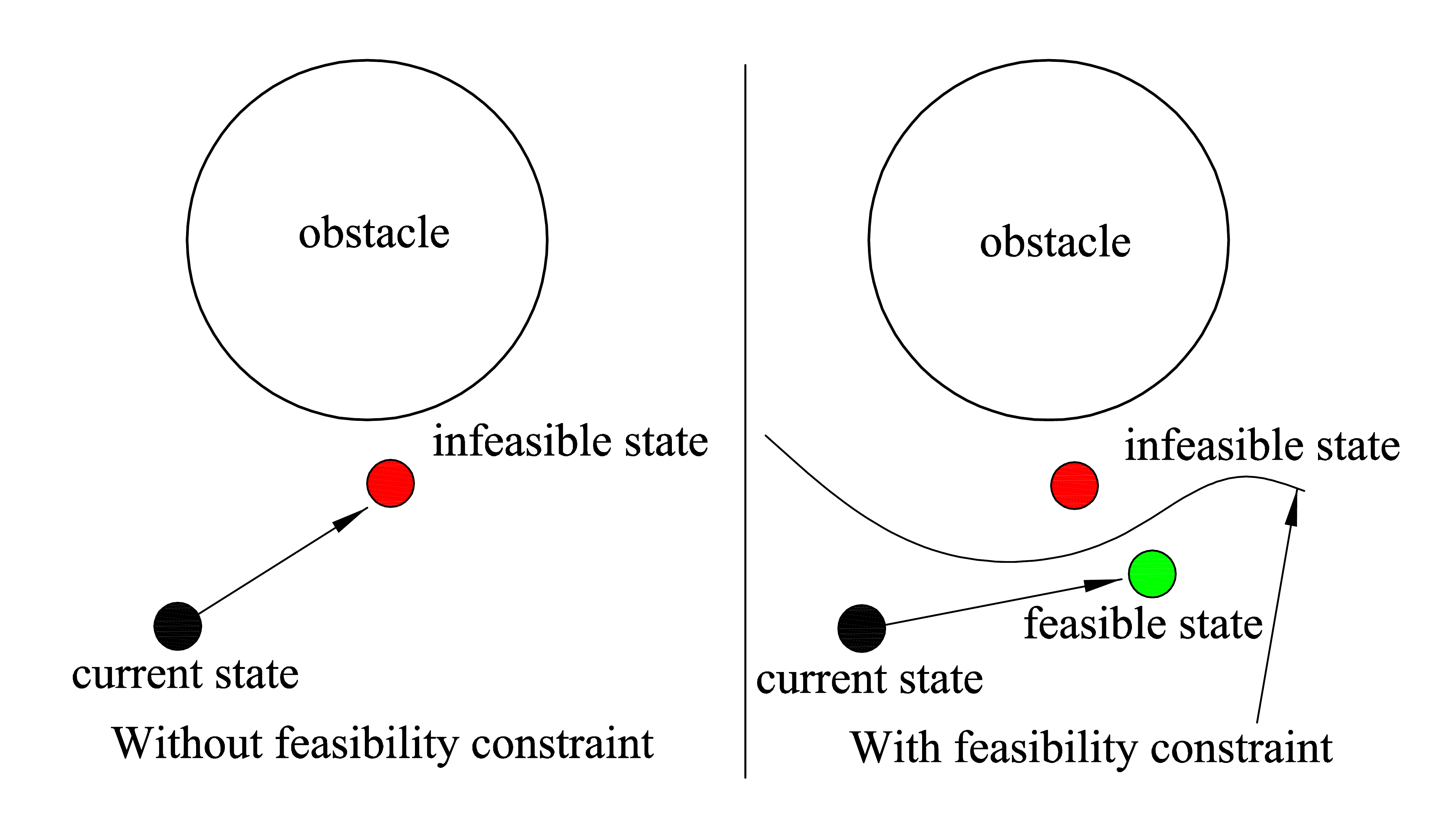}
	\caption{An illustration of how a feasibility constraint works for a robot control problem. A feasibility constraint prevents the robot from going into the infeasible state.}	
	\label{fig:illustration}
\end{figure}

An illustrative example of how a feasibility constraint works is shown in Fig. \ref{fig:illustration}. A robot whose control is determined by solving the QP (\ref{eqn:prob_qp}), subject to (\ref{eqn:prob_cbf}), (\ref{eqn:prob_clf}) and (\ref{eqn:prob_ctrl}), will run close to an obstacle in the following step. The next state may be infeasible for the QP associated with that next step. For example, the state denoted by the red dot in Fig. \ref{fig:illustration} may have large speed such that the robot cannot find a control to avoid the obstacle in the next step. If a feasibility constraint can prevent the robot from reaching this state, then the QP is feasible.

After we find a feasibility constraint, we can enforce it through a CBF and take it as an additional constraint for (\ref{eqn:prob_qp}) to guarantee the feasibility given system state $\bm x$. We show how we can determine an appropriate feasibility constraint in the following section.

\section{Feasibility Constraint}

\label{sec:fc}

We begin with a simple example to illustrate the necessity for a feasibility constraint for the CBF-CLF based QPs.

\subsection{Example: Adaptive Cruise Control}

\label{sec:hocbf:moti}

Consider the adaptive cruise control (ACC) problem with the ego
(controlled) vehicle dynamics in the form: 
\begin{equation}
\underbrace{\left[
	\begin{array}
	[c]{c}%
	\dot{v}(t)\\
	\dot{z}(t)
	\end{array}
	\right]  }_{\dot{\bm x}(t)}=\underbrace{\left[
	\begin{array}
	[c]{c}%
	-\frac{1}{M}F_{r}(v(t))\\
	v_p - v(t)
	\end{array}
	\right]  }_{f(\bm x(t))}+\underbrace{\left[
	\begin{array}
	[c]{c}%
	\frac{1}{M}\\
	0
	\end{array}
	\right]  }_{g(\bm x(t))}u(t)\label{eqn:veh}%
\end{equation}
where $M$ denotes the mass of the ego vehicle, $z(t)$ denotes the distance between the preceding and the ego vehicles, $v_p \geq 0, v(t) \geq 0$ denote the speeds of the preceding and the ego vehicles, respectively, and $F_{r}(v(t))$
denotes the resistance force, which is expressed \cite{Khalil2002} as:
$$F_{r}(v(t))=f_{0}sgn(v(t))+f_{1}v(t)+f_{2}v^{2}(t),$$ where $f_{0}>0,$
$f_{1}>0$ and $f_{2}>0$ are scalars determined empirically. The first term in
$F_{r}(v(t))$ denotes the Coulomb friction force, the second term denotes the
viscous friction force and the last term denotes the aerodynamic drag. The control $u(t)$ is the driving force of the ego vehicle subject to the
constraint:
\begin{equation}
-c_dMg\leq u(t)\leq c_aMg,\forall t\geq0, \label{eqn:scons}%
\end{equation}
where $c_a>0$ and $c_d>0$ are the maximum acceleration and deceleration coefficients, respectively, and $g$ is the gravity constant.

We require that the distance $z(t)$ between the ego vehicle and its immediately
preceding vehicle be greater than $l_{0} > 0$,
i.e.,
\begin{equation}
\label{eqn:safety}z(t) \geq l_{0},\forall t\geq0.
\end{equation}

Let $b(\bm x(t)):=z(t)-l_{0}$. The relative
degree of $b(\bm x(t))$ is $m=2$, so we choose a HOCBF following Def.
\ref{def:hocbf} by defining $\psi_{0}(\bm x(t)):=b(\bm x(t))$, $\alpha
_{1}(\psi_{0}(\bm x(t))):=p_1\psi_{0}(\bm x(t))$ and $\alpha_{2}(\psi_{1}(\bm
x(t))):=p_2\psi_{1}(\bm x(t)), p_1 > 0, p_2 > 0$. We then seek a control for the ego vehicle such
that the constraint (\ref{eqn:safety}) is satisfied. The control $u(t)$ should
satisfy (\ref{eqn:constraint}) which in this case is: 
\begin{equation}
\begin{aligned} \underbrace{\frac{F_r(v(t))}{M}}_{L_f^2b(\bm x(t))} + \underbrace{\frac{-1}{M}}_{L_gL_fb(\bm x(t))}\times u(t) + \underbrace{p_1(v_p - v(t))}_{S(b(\bm x(t)))} \\+ \underbrace{p_2(v_p - v(t)) + p_1p_2(z(t) - l_0)}_{\alpha_2(\psi_1(\bm x(t)))} \geq 0. \end{aligned}
\label{eqn:safety_ex2}%
\end{equation}

Suppose we wish to minimize $\int_{0}^{T}\left(\frac{u(t)-F_{r}(v(t))}{M}\right)^{2}dt$, in which case we have a constrained optimal control problem. We can then use the
QP-based method introduced at the end of the last section to solve this ACC
problem. However, the HOCBF constraint (\ref{eqn:safety_ex2}) can easily
conflict with $-c_dMg\leq u(t)$ in (\ref{eqn:scons}), i.e., the ego vehicle cannot brake in time under control constraint (\ref{eqn:control}) so that the safety constraint (\ref{eqn:safety}) is satisfied when the two vehicles get close to each other. This is intuitive when we rewrite (\ref{eqn:safety_ex2}) in the form:
\begin{equation}\label{eqn:safety_ex22}
\begin{aligned}  \frac{1}{M}u(t)\!\leq \!\frac{F_r(v(t))}{M} \!+ \!(p_1\! +\! p_2)(v_p \!-\! v(t)) \!+\! p_1p_2(z(t) \!-\! l_0). \end{aligned}
\end{equation}
The right-hand side above is usually negative when the two vehicles get close to each other. If it is smaller than $-c_dMg$, the HOCBF constraint (\ref{eqn:safety_ex2}) will conflict with $-c_dMg\leq u(t)$ in (\ref{eqn:scons}). When this happens, the QP will be infeasible. In the rest of the paper, we show how we can solve this infeasibility problem in general by a feasibility constraint as in Def. \ref{def:fc}.

\subsection{Feasibility Constraint for Relative-Degree-One Safety Constraints}

For simplicity, we start with feasibility constraints for a relative-degree-one safety constraint.

Suppose we have a constraint $b(\bm x)\geq 0$ with relative degree one for system (\ref{eqn:affine}), where $b:\mathbb{R}^n\rightarrow \mathbb{R}$. Then we can define $b(\bm x)$ as a  HOCBF with $m = 1$ as in Def. \ref{def:hocbf}, i.e., we have a ``traditional'' CBF. Following (\ref{eqn:constraint}), any control $\bm u\in U$ should satisfy the CBF constraint:
\begin{equation} \label{eqn:cbf}
- L_gb(\bm x)\bm u \leq L_fb(\bm x)+ \alpha(b(\bm x)),
\end{equation}
where $\alpha(\cdot)$ is a class $\mathcal{K}$ function of its argument. We define a set of controls that satisfy the last equation as:
\begin{equation}\label{eqn:cbfset}
K(\bm x) = \{\bm u\in \mathbb{R}^q:-L_gb(\bm x)\bm u \leq L_fb(\bm x) + \alpha(b(\bm x))\}.
\end{equation}

Our analysis for determining a feasibility constraint depends on whether any component of the vector $L_gb(\bm x)$ will change sign in the time interval $[0,T]$ or not.

\subsubsection{All components in $L_gb(\bm x)$ do not change sign} Since all components in $L_gb(\bm x)$ do not change sign for all $\bm x\in X$, the inequality constraint for each control component does not change sign if we multiply each component of $L_gb(\bm x)$ by the corresponding one of the control bounds in (\ref{eqn:control}). Therefore, we assume that $L_gb(\bm x)\leq \bm 0 \text{ (componentwise)}, \bm 0\in\mathbb{R}^q$  in the rest of this section. The analysis for other cases (each component of $L_gb(\bm x)$ is either non-negative or non-positive) is similar. Not all the components in $L_gb(\bm x)$ can be 0 due to the relative degree definition in Def. \ref{def:relative}. We can multiply  the control bounds (\ref{eqn:control}) by the vector $-L_gb(\bm x)$, and get
\begin{equation}\label{eqn:ctrl_trans}
- L_gb(\bm x)\bm u_{\min}\leq - L_gb(\bm x)\bm u \leq - L_gb(\bm x)\bm u_{\max},
\end{equation}
The control constraint (\ref{eqn:ctrl_trans}) is actually a relaxation of the control bound (\ref{eqn:control}) as we multiply each component of $L_gb(\bm x)$ by the corresponding one of the control bounds in (\ref{eqn:control}), and then add them together. We define 
\begin{equation} \label{eqn:set_ex}
\begin{aligned}U_{ex}(\bm x) &=\{\bm u\in\mathbb{R}^q:\\&- L_gb(\bm x)\bm u_{\min}\leq - L_gb(\bm x)\bm u \leq - L_gb(\bm x)\bm u_{\max}\},\end{aligned} 
\end{equation}
It is obvious that $U$ is a subset of $U_{ex}(\bm x)$. Nonetheless, the relaxation set $U_{ex}(\bm x)$ does not negatively affect the property of the following lemma:

\begin{lemma}\label{lem:equiv}
	If the control $\bm u$ is such that (\ref{eqn:ctrl_trans}) is conflict-free with (\ref{eqn:cbf}) for all $\bm x\in X$, then the control bound (\ref{eqn:control}) is also conflict-free with (\ref{eqn:cbf}).
\end{lemma}

\textbf{Proof:} Let $g = (g_1,\dots,g_q)$ in (\ref{eqn:affine}), where $g_i : \mathbb{R}^n\rightarrow \mathbb{R}^n, i,\in\{1,\dots,q\}$. We have that $L_gb(\bm x) = (L_{g_1}b(\bm x),\dots, L_{g_q}b(\bm x)) \in\mathbb{R}^{1\times q}$. For the control bound $u_{i,min}\leq u_i\leq u_{i,max}, i\in\{1,\dots,q\}$ in (\ref{eqn:control}), we can multiply by $-L_{g_i}b(\bm x)$ and get 
$$\begin{aligned}-L_{g_i}b(\bm x)u_{i,min}\leq -L_{g_i}b(\bm x)u_i\leq -L_{g_i}b(\bm x)u_{i,max},\\ i\in\{1,\dots,q\},\end{aligned}$$
as we have assumed that $L_gb(\bm x)\leq \bm 0$. If we take the summation of the inequality above over all $i\in\{1,\dots, q\}$, then we obtain the constraint (\ref{eqn:ctrl_trans}). Therefore, the satisfaction of (\ref{eqn:control}) implies the satisfaction of (\ref{eqn:ctrl_trans}). 
Then $U$ defined in (\ref{eqn:control}) is a subset of $U_{ex}(\bm x)$. It is obvious that the boundaries of the set $U_{ex}(\bm x)$ in (\ref{eqn:set_ex}) and $K(\bm x)$ in (\ref{eqn:cbfset}) are hyperplanes, and these boundaries are parallel to each other for all $\bm x\in X$. Meanwhile, the two boundaries of $U_{ex}(\bm x)$ pass through the two corners $\bm u_{\min}, \bm u_{\max}$ of the set $U$ (a polyhedron) following (\ref{eqn:set_ex}), respectively.  If there exists a control $\bm u_1 \in U_{ex}(\bm x)$ that satisfies (\ref{eqn:cbf}), then the boundary of the set $K(\bm x)$ in (\ref{eqn:cbfset}) lies either between the two hyperplanes defined by $U_{ex}(\bm x)$ or above these two hyperplanes (i.e., $U_{ex}(\bm x)$ is a subset of $K(\bm x)$ in (\ref{eqn:cbfset})). In the latter case, this lemma is true as $U$ is a subset of $U_{ex}(\bm x)$. In the former case, we can always find another control $\bm u_2\in U$ that satisfies (\ref{eqn:cbf}) as the boundary of $K(\bm x)$ in (\ref{eqn:cbfset}) is parallel to the two $U_{ex}(\bm x)$ boundaries that respectively pass through the two corners $\bm u_{\min}, \bm u_{\max}$ of the set $U$. Therefore, although $U$ is a subset of $U_{ex}(\bm x)$, it follows that if (\ref{eqn:ctrl_trans}) is conflict-free with (\ref{eqn:cbf}) in terms of $\bm u$ for all $\bm x\in X$, the control bound (\ref{eqn:control}) is also conflict-free with (\ref{eqn:cbf}). $\qquad\qquad\qquad\qquad\qquad\qquad\qquad\blacksquare$

As motivated by Lem. \ref{lem:equiv}, in order to determine if (\ref{eqn:cbf}) complies with (\ref{eqn:control}), we may just consider (\ref{eqn:cbf}) and (\ref{eqn:ctrl_trans}). Since there are two inequalites in (\ref{eqn:ctrl_trans}), we have two cases to consider: $(i) - L_gb(\bm x)\bm u \leq - L_gb(\bm x)\bm u_{\max}$ and (\ref{eqn:cbf}); $(ii) - L_gb(\bm x)\bm u_{\min}\leq - L_gb(\bm x)\bm u$ and (\ref{eqn:cbf}). It is obvious that there always exists a control $\bm u$ such that the two inequalities in case $(i)$ are satisfied for all $\bm x\in X$, while this may not be true for case $(ii)$, depending on $\bm x$. Therefore, in terms of avoiding the conflict between the CBF constraint (\ref{eqn:cbf}) and (\ref{eqn:ctrl_trans}) that leads to the infeasibility of problem (\ref{eqn:prob_qp}), subject to (\ref{eqn:prob_cbf})-(\ref{eqn:prob_ctrl}), we wish to satisfy:
\begin{equation}\label{eqn:fea_r1}
L_fb(\bm x) + \alpha(b(\bm x))\geq -L_gb(\bm x)\bm u_{\min}.
\end{equation}
This is called the {\bf feasibility constraint} for problem (\ref{eqn:prob_qp}), subject to (\ref{eqn:prob_cbf})-(\ref{eqn:prob_ctrl}) in the case of a relative-degree-one safety constraint $b(\bm x)\geq 0$ in (\ref{eqn:safetycons}).

The relative degree of the feasibility constraint (\ref{eqn:fea_r1}) is also one with respect to dynamics (\ref{eqn:affine}) as we have $b(\bm x)$ in it. In order to find a control such that the feasibility constraint (\ref{eqn:fea_r1}) is guaranteed to be satisfied, we define
\begin{equation}\label{eqn:bF}
b_F(\bm x) = L_fb(\bm x) + \alpha(b(\bm x)) + L_gb(\bm x)\bm u_{\min} \geq 0,
\end{equation}
so that $b_F(\bm x)$ is a CBF as in Def. \ref{def:hocbf}.
Then, we can get a feedback controller $K_F(\bm x)$ that guarantees the CBF constraint (\ref{eqn:cbf}) and the control bounds (\ref{eqn:control}) do not conflict with each other:
\begin{equation}\label{eqn:cbf_fea}
K_{F}(\bm x) = \{\bm u\in \mathbb{R}^q: L_fb_F(\bm x) + L_gb_F(\bm x)\bm u + \alpha_f(b_F(\bm x))\geq 0\},
\end{equation}
if $b_F(\bm x(0))\geq 0$, where $\alpha_f(\cdot)$ is a class $\mathcal{K}$ function.

\begin{theorem}\label{thm:fea}
	If Problem \ref{prob:general} is initially feasible and the CBF constraint in (\ref{eqn:cbf_fea}) corresponding to (\ref{eqn:fea_r1}) does not conflict with both the control bounds (\ref{eqn:control}) and (\ref{eqn:cbf}) at the same time, any controller $\bm u \in K_F(\bm x)$ guarantees the feasibility of problem (\ref{eqn:prob_qp}), subject to (\ref{eqn:prob_cbf})-(\ref{eqn:prob_ctrl}).
\end{theorem}
\textbf{Proof:} If Problem \ref{prob:general} is initially feasible, then the CBF constraint (\ref{eqn:cbf}) for the safety requirement (\ref{eqn:safetycons}) does not conflict with the control bounds (\ref{eqn:control}) at time 0. It also does not conflict with the constraint (\ref{eqn:ctrl_trans}) as $U$ is a subset of $U_{ex}(\bm x)$ that is defined in (\ref{eqn:set_ex}). In other words, $b_F(\bm x(0))\geq 0$ holds in the feasibility constraint (\ref{eqn:fea_r1}). Thus, the initial condition for the CBF in Def. \ref{def:hocbf} is satisfied. By Thm. \ref{thm:hocbf}, we have that  $b_F(\bm x(t))\geq 0, \forall t\geq 0$. Therefore, the CBF constraint (\ref{eqn:cbf}) does not conflict with the constraint (\ref{eqn:ctrl_trans}) for all $t\geq 0$. By Lemma \ref{lem:equiv}, the CBF constraint (\ref{eqn:cbf}) also does not conflict with the control bound (\ref{eqn:control}). Finally, since the CBF constraint in (\ref{eqn:cbf_fea}) corresponding to (\ref{eqn:fea_r1}) does not conflict with the control bounds (\ref{eqn:control}) and (\ref{eqn:cbf})  at the same time by assumption, we conclude that the feasibility of the problem is guaranteed. $\blacksquare$

The condition ``the CBF constraint in (\ref{eqn:cbf_fea}) corresponding to (\ref{eqn:fea_r1}) does not conflict with both the control bounds (\ref{eqn:control}) and (\ref{eqn:cbf})  at the same time'' in Thm. \ref{thm:fea} is too strong. If this condition is not satisfied, then the problem can still be infeasible. In order to relax this condition, one option is to recursively define other new feasibility constraints for the feasibility constraint (\ref{eqn:fea_r1}) to address the possible conflict between (\ref{eqn:cbf_fea}) and (\ref{eqn:control}), and (\ref{eqn:cbf}). However, the number of iterations is not bounded, and we may have a large (unbounded) set of feasibility constraints. 

In order to address the unbounded iteration issue in finding feasibility constraints, we can try to express the feasibility constraint in (\ref{eqn:cbf_fea}) so that it is in a form which is similar to that of the CBF constraint (\ref{eqn:cbf}). If this is achieved, we can make these two constraints compliant with each other, and thus address the unbounded iteration issue mentioned above. Therefore, we try to construct the CBF constraint in (\ref{eqn:cbf_fea}) so that it takes the form:
\begin{equation} \label{eqn:cbf_fea_re}
L_fb(\bm x) + L_gb(\bm x)\bm u + \alpha(b(\bm x)) + \varphi(\bm x, \bm u) \geq 0
\end{equation}
for some appropriately selected function $\varphi(\bm x, \bm u)$. One obvious choice for $\varphi(\bm x, \bm u)$ immediately following (\ref{eqn:cbf_fea}) is $\varphi(\bm x, \bm u) = L_fb_F(\bm x) + L_gb_F(\bm x)\bm u + \alpha_f(b_F(\bm x)) - L_fb(\bm x) - L_gb(\bm x)\bm u - \alpha(b(\bm x))$, which can be simplified through a proper choice of the class $\mathcal{K}$ functions $\alpha(\cdot), \alpha_f(\cdot)$, as will be shown next. Since we will eventually include the constraint $\varphi(\bm x, \bm u)\geq 0$ into our QPs (shown later) to address the infeasibility problem, we wish its relative degree to be low. Otherwise, it becomes necessary to use HOCBFs to make the control show up in enforcing $\varphi(\bm x)\geq 0$ (instead of $\varphi(\bm x, \bm u)\geq 0$ due to its high relative degree), which could make the corresponding HOCBF constraint complicated, and make it easily conflict with the control bound (\ref{eqn:control}) and the CBF constraint (\ref{eqn:cbf}), and thus leading to the infeasibility of the QPs.  Therefore,  we define a candidate function as follows (note that a relative-degree-zero function means that the control $\bm u$ directly shows up in the function itself):

\begin{definition} [Candidate $\varphi(\bm x, \bm u)$ function]
	A function $\varphi(\bm x, \bm u)$ in (\ref{eqn:cbf_fea_re}) is a {\it candidate function} if its relative degree with respect to (\ref{eqn:affine}) is either one or zero. 
\end{definition}
\textbf{Finding candidate $\varphi(\bm x, \bm u)$:}
	In order to find a candidate $\varphi(\bm x, \bm u)$ from the reformulation of the CBF constraint in (\ref{eqn:cbf_fea}), we can properly choose the class $\mathcal{K}$ function $\alpha(\cdot)$ in (\ref{eqn:cbf}). A typical choice for $\alpha(\cdot)$ is a linear function, in which case we automatically have the constraint formulation (\ref{eqn:cbf_fea_re}) by substituting the function $b_F(\bm x)$ from (\ref{eqn:bF}) into (\ref{eqn:cbf_fea}), and get 
	$$
	\begin{aligned}
	\varphi(\bm x, \bm u) = L_f^2b(\bm x) + L_gL_fb(\bm x) \bm u + L_f(L_gb(\bm x)\bm u_{min}) \\+ L_g(L_gb(\bm x)\bm u_{min})\bm u + \alpha_f(b_F(\bm x))- b(\bm x).
	\end{aligned}$$ 
	Note that it is possible that $L_gL_fb(\bm x) = 0$ and $L_g(L_gb(\bm x)\bm u_{min}) = 0$ (depending on the dynamics (\ref{eqn:affine}) and the CBF $b(\bm x)$), in which case the relative degree of $\varphi(\bm x, \bm u)$ (written as $\varphi(\bm x)$) is one as we have $\alpha_f(b_F(\bm x))$ in it and $b_F(\bm x)$ is a function of $b(\bm x)$. 

If the relative degree of $\varphi(\bm x, \bm u)$ is zero (e.g., $L_gL_fb(\bm x) = 0$ and $L_g(L_gb(\bm x)\bm u_{min}) = 0$ are not satisfied above), we wish to require that 
\begin{equation} \label{eqn:r0}
\varphi(\bm x, \bm u)\geq 0, 
\end{equation}
such that the satisfaction of the CBF constraint (\ref{eqn:cbf}) implies the satisfaction of the CBF constraint (\ref{eqn:cbf_fea_re}), and the satisfaction of the CBF constraint (\ref{eqn:cbf_fea_re}) implies the satisfaction of (\ref{eqn:fea_r1}) by Thm. \ref{thm:hocbf}, i.e., the CBF constraint (\ref{eqn:cbf}) does not conflict with the control bound (\ref{eqn:control}). Besides, if (\ref{eqn:r0}) happens to not conflict with both (\ref{eqn:cbf}) and (\ref{eqn:control}) at the same time, depending on the CBF $b(\bm x)$ and the dynamics (\ref{eqn:affine}), then the QPs are guaranteed to be feasible. The constraint (\ref{eqn:r0}) is simplier than (\ref{eqn:cbf_fea}) as all the terms in the CBF constraint (\ref{eqn:cbf}) are removed through (\ref{eqn:cbf_fea_re}), thus, it is less likely to conflict with the CBF constraint (\ref{eqn:cbf}) and the control bound (\ref{eqn:control}) in the QP. This is more helpful in the case of safety constraints with high relative degree (in the next subsection) as the HOCBF constraint (\ref{eqn:constraint}) has many complicated terms, and it is better to remove these terms in the feasibility constraint and just consider (\ref{eqn:r0}) in the QP in order to make (\ref{eqn:r0}) compliant with (\ref{eqn:cbf}) and (\ref{eqn:control}).

If the relative degree of a candidate $\varphi(\bm x, \bm u)$ with respect to (\ref{eqn:affine}) is one, i.e., $\varphi(\bm x, \bm u) \equiv \varphi(\bm x)$, we define a set $U_{s}(\bm x)$:
\begin{equation}\label{eqn:set_single}
U_{s}(\bm x) = \{\bm u\in\mathbb{R}^q: L_f\varphi(\bm x) + L_g\varphi(\bm x)\bm u + \alpha_u(\varphi(\bm x))\geq 0\}.
\end{equation}
where $\alpha_u(\cdot)$ is a class $\mathcal{K}$ function. 

From the set of candidate functions $\varphi(\bm x)$, if we can find one that satisfies the conditions of the following theorem, then the feasibility of problem (\ref{eqn:prob_qp}), subject to (\ref{eqn:prob_cbf})-(\ref{eqn:prob_ctrl}) is guaranteed:
\begin{theorem} \label{thm:feasible_unique}
	If $\varphi(\bm x)$ is a candidate function such that $\varphi(\bm x(0))\geq 0, L_f\varphi(\bm x) \geq 0$, $L_g\varphi(\bm x) = \gamma L_gb(\bm x)$, for some $\gamma > 0, \forall \bm x\in X$ and $\bm 0\in U$, then any controller $\bm u(t) \in U_s(\bm x), \forall t\geq 0$ guarantees the feasibility of problem (\ref{eqn:prob_qp}), subject to (\ref{eqn:prob_cbf})-(\ref{eqn:prob_ctrl}).
\end{theorem}
\textbf{Proof:} Since $\varphi(\bm x)$ is a candidate function, we can define a set $U_s(\bm x)$ as in (\ref{eqn:set_single}). If $\varphi(\bm x(0))\geq 0$ and $\bm u(t) \in U_s(\bm x), \forall t\geq 0$, we have that $\varphi(\bm x(t))\geq 0, \forall t\geq 0$ by Thm. \ref{thm:hocbf}. Then, the satisfaction of the CBF constraint (\ref{eqn:cbf}) corresponding to the safety constraint (\ref{eqn:safetycons}) implies the satisfaction of the CBF constraint (\ref{eqn:cbf_fea_re}) (equivalent to (\ref{eqn:cbf_fea})) for the feasibility constraint (\ref{eqn:fea_r1}). In other words, the CBF constraint (\ref{eqn:cbf}) automatically guarantees that it will not conflict with the control constraint (\ref{eqn:ctrl_trans}) as the satisfaction of (\ref{eqn:cbf_fea_re}) implies the satisfaction of (\ref{eqn:fea_r1}) following Thm. \ref{thm:hocbf} and (\ref{eqn:fea_r1}) guarantees that (\ref{eqn:cbf}) and (\ref{eqn:ctrl_trans}) are conflict-free. By Lemma \ref{lem:equiv}, the CBF constraint (\ref{eqn:cbf}) will also not conflict with the control bound $U$ in (\ref{eqn:control}), i.e. $K(\bm x)\cap U\ne \emptyset$, where $K(\bm x)$ is defined in (\ref{eqn:cbfset}).

Since $L_f\varphi(\bm x) \geq 0$, we have that $\bm 0 \in U_s(\bm x)$. We also have $\bm 0 \in U(\bm x)$, thus, $U_s(\bm x) \cap U \ne \emptyset$ is guaranteed. Since $L_g\varphi(\bm x) = \gamma L_gb(\bm x), \gamma > 0$, the two hyperplanes of the two half spaces formed by $U_s(\bm x)$ in (\ref{eqn:set_single}) and $K(\bm x)$ in (\ref{eqn:cbfset}) are parallel to each other, and the normal directions of the two hyperplanes along the half space direction are the same. Thus, $U_s(\bm x)\cap K(\bm x)$ is either $U_s(\bm x)$ or $K(\bm x)$, i.e., $U_s(\bm x)\cap K(\bm x)\cap U$ equals either $U_s(\bm x) \cap U$ or $K(\bm x) \cap U$. As $U_s(\bm x) \cap U \ne \emptyset$ and $K(\bm x)\cap U\ne \emptyset$, we have $U_s(\bm x)\cap K(\bm x)\cap U \ne \emptyset, \forall \bm x\in X$. Therefore, the CBF constraint (\ref{eqn:cbf}) does not conflict with the control bound (\ref{eqn:control}) and the CBF constraint in $U_s(\bm x)$ at the same time, and we can conclude that the problem is guaranteed to be feasible. $\qquad\qquad\qquad\qquad\;\;\;\blacksquare$

	The conditions in Thm. \ref{thm:feasible_unique} are {\bf sufficient conditions} for the feasibility of problem (\ref{eqn:prob_qp}), subject to (\ref{eqn:prob_cbf})-(\ref{eqn:prob_ctrl}). Under the conditions in Thm \ref{thm:feasible_unique}, we can claim that $\varphi(\bm x) \geq 0$ is a single {\bf feasibility constraint} that guarantees the feasibility of problem (\ref{eqn:prob_qp}), subject to (\ref{eqn:prob_cbf})-(\ref{eqn:prob_ctrl}) in the case that the safety constraint (\ref{eqn:safetycons}) is with relative degree one (i.e., $m = 1$ in (\ref{eqn:prob_cbf})). 

\textbf{Finding valid $\varphi(\bm x)$:}
		A valid $\varphi(\bm x)$ is a function that satisfies the conditions in Thm. \ref{thm:feasible_unique}. The conditions in Thm. \ref{thm:feasible_unique} may be conservative, and how to determine such a $\varphi(\bm x)$ function is the remaining problem. For a general system (\ref{eqn:affine}) and safety constraint (\ref{eqn:safetycons}), we can parameterize the definition of the CBF (\ref{eqn:cbf}) for the safety and the CBF constraint for the feasibility constraint (\ref{eqn:cbf_fea}), i.e., parameterize $\alpha(\cdot)$ and $\alpha_F(\cdot)$, such as the form in \cite{Xiao2020CDC}, and then choose the parameters to satisfy the conditions in Thm. \ref{thm:feasible_unique}. 
		
		\begin{remark}An example for determining such a $\varphi(\bm x)$ for the ACC problem in Sec. \ref{sec:hocbf:moti} can be found in the end of this section.  However, it is still not guaranteed that such $\varphi(\bm x)$ functions can be found. To address this, we may consider a special class of dynamics (\ref{eqn:affine}), and then formulate a systematic way to derive such $\varphi(\bm x)$ functions. In the case of such dynamics, we may even relax some of the conditions in Thm. \ref{thm:feasible_unique}. For example, if $g(\bm x)$ in (\ref{eqn:affine}) is independent of $\bm x$ and the safety constraint (\ref{eqn:safetycons}) is in linear form, then it is very likely that the condition $L_g\varphi(\bm x) = \gamma L_gb(\bm x)$, for some $\gamma > 0$ in Thm. \ref{thm:feasible_unique} is satisfied, and thus this condition may be removed.  
        \end{remark}

We can now get a feasible problem from the original problem (\ref{eqn:prob_qp}), subject to (\ref{eqn:prob_cbf})-(\ref{eqn:prob_ctrl}) in the form:
\begin{equation}\label{eqn:prob_qp1}
\min_{\bm u(t), \delta(t)} \int_{0}^T ||\bm u(t)||^2 + p\delta^2(t) dt
\end{equation}
subject to
the feasibility constraint (\ref{eqn:r0}) if the relative degree of $\varphi(\bm x, \bm u)$ is 0; otherwise, subject to the CBF constraint in (\ref{eqn:set_single}). The cost (\ref{eqn:prob_qp1}) is also subject to the CBF constraint (\ref{eqn:cbf}), the control bound (\ref{eqn:control}), and the CLF constraint:
\begin{equation} \label{eqn:prob_clf1}
L_{f}V(\bm x)+L_{g}V(\bm x)\bm u + \epsilon V(\bm x)\leq \delta(t),
\end{equation}
where $\varphi(\bm x)$ satisfies the conditions in Thm. \ref{thm:feasible_unique} for (\ref{eqn:set_single}), and (\ref{eqn:r0}) is assumed to be non-conflicting with the CBF constraint (\ref{eqn:cbf}) and  the control bound (\ref{eqn:control}) at the same time. In order to guarantee feasibility, we may try to find a $\varphi(\bm x)$ that has relative degree one, and that satisfies the conditions in Thm. \ref{thm:feasible_unique}.

\subsubsection{Some Components in $L_gb(\bm x)$ Change Sign} Recall that $L_gb(\bm x) = (L_{g_1}b(\bm x), \dots, L_{g_q}b(\bm x))\in\mathbb{R}^{1\times q}$. If $L_{g_i}b(\bm x), i\in\{1,\dots,q\}$ changes sign in $[0,T]$, then we have the following symmetric and non-symmetric cases to consider in order to find a valid feasibility constraint.

Let $\bm u = (u_1,\dots, u_q)$, $\bm u_{\min} = (u_{1,\min},\dots, u_{q,\min}) \leq \bm 0$, $\bm u_{\max} = (u_{1,\max},\dots, u_{q,\max})\geq\bm 0, \bm 0\in\mathbb{R}^q$. 

\textbf{Case 1:} the control bound for $u_{i}, i\in\{1,\dots,q\}$ is symmetric, i.e. $u_{i,\max} = -u_{i,\min}$. In this case, by multiplying $-L_{g_i}b(\bm x)$ by the control bound for $u_i$, we have
\begin{equation} \label{eqn:ctrl_single}
-L_{g_i}b(\bm x)u_{i,\min} \leq -L_{g_i}b(\bm x)u_i\leq -L_{g_i}b(\bm x)u_{i,\max}
\end{equation}
if $L_{g_i}b(\bm x) < 0$. When $L_{g_i}b(\bm x)$ changes sign at some time $t_1\in[0,T]$, then the sign of the last equation will be reversed. However, since $u_{i,\max} = -u_{i,\min}$, we have exactly the same constraint as (\ref{eqn:ctrl_single}), and $-L_{g_i}b(\bm x)u_{i,\min}$ will still be continuously differentiable when we construct the feasibility constraint as in (\ref{eqn:fea_r1}). Therefore, the feasibility constraint (\ref{eqn:fea_r1}) will not be affected by the sign change of $L_{g_i}b(\bm x), i\in\{1,\dots,q\}$.

\textbf{Case 2:} the control bound for $u_{i}, i\in\{1,\dots,q\}$ is not symmetric, i.e., $u_{i,\max} \ne -u_{i,\min}$. In this case, we can define:
\begin{equation}\label{eqn:ctrl_lim}
u_{i,\lim} := \min\{|u_{i,\min}|, u_{i,\max}\}
\end{equation}

Considering (\ref{eqn:ctrl_lim}), we have the following constraint
\begin{equation}
-u_{i,\lim} \leq u_i\leq u_{i,\lim}.
\end{equation}
The satisfaction of the last equation implies the satisfaction of $u_{i,\min}\leq u_i\leq u_{i,\max}$ in (\ref{eqn:control}).

If $L_{g_i}b(\bm x)< 0$, we multiply the control bound by $-L_{g_i}b(\bm x)$ for $u_i$ and have the following constraint
\begin{equation} \label{eqn:ctrl_single2}
L_{g_i}b(\bm x)u_{i,\lim} \leq -L_{g_i}b(\bm x)u_i\leq -L_{g_i}b(\bm x)u_{i,\lim}
\end{equation}
The satisfaction of (\ref{eqn:ctrl_single2}) implies the satisfaction of (\ref{eqn:ctrl_single}) following (\ref{eqn:ctrl_lim}). Now, the control bound for $u_i$ is converted to the symmetric case, and the feasibility constraint (\ref{eqn:fea_r1}) will not be affected by the sign change of $L_{g_i}b(\bm x), i\in\{1,\dots,q\}$.

\subsection{Feasibility Constraint for High-Relative-Degree Safety Constraints}

Suppose we have a constraint $b(\bm x)\geq 0$ with relative degree $m\geq 1$ for system (\ref{eqn:affine}), where $b:\mathbb{R}^n\rightarrow \mathbb{R}$. Then we can define $b(\bm x)$ as a HOCBF as in Def. \ref{def:hocbf}. Any control $\bm u\in U$ should satisfy the HOCBF constraint (\ref{eqn:constraint}).

In this section, we also assume that $L_gL_f^{m-1}b(\bm x) \leq \bm 0, \bm 0\in\mathbb{R}^q$ and all components in $L_gL_f^{m-1}b(\bm x)$ do not change sign in $[0,T]$. The analysis for all other cases is similar to the last subsection.

Similar to (\ref{eqn:cbf}), we rewrite the HOCBF constraint (\ref{eqn:constraint}) as
\begin{equation} \label{eqn:hocbf}
-L_gL_f^{m-1}b(\bm x)\bm u\leq L_f^{m}b(\bm x) +  S(b(\bm x)) + \alpha_m(\psi_{m-1}(\bm x))
\end{equation}

We can multiply the control bounds (\ref{eqn:control}) by the vector $-L_gL_f^{m-1}b(\bm x)$:
\begin{equation}\label{eqn:ctrl_transH}
\begin{aligned}
- L_gL_f^{m-1}b(\bm x)\bm u_{\min}\leq - L_gL_f^{m-1}b(\bm x)\bm u \\\leq - L_gL_f^{m-1}b(\bm x)\bm u_{\max},
\end{aligned}
\end{equation}
As in (\ref{eqn:ctrl_trans}), the last equation is also a relaxation of the original control bound (\ref{eqn:control}), and Lem. \ref{lem:equiv} still applies in the high-relative-degree-constraint case.

The HOCBF constraint (\ref{eqn:hocbf}) may conflict with the left inequality of the transformed control bound (\ref{eqn:ctrl_transH}) when its right hand side is smaller than $- L_gL_f^{m-1}b(\bm x)\bm u_{\min}$. Therefore, we wish to have
\begin{equation} \label{eqn:feaH}
\begin{aligned}
L_f^{m}b(\bm x) +  S(b(\bm x)) + \alpha_m(\psi_{m-1}(\bm x)) \geq -L_gL_f^{m-1}b(\bm x)\bm u_{\min}.
\end{aligned}
\end{equation}
This is called the {\bf feasibility constraint} for the problem (\ref{eqn:prob_qp}), subject to (\ref{eqn:prob_cbf})-(\ref{eqn:prob_ctrl}) in the case of a high-relative-degree constraint $b(\bm x)\geq 0$ in (\ref{eqn:safetycons}).

In order to find a control such that the feasibility constraint (\ref{eqn:fea_r1}) is guaranteed to be satisfied, we define 
$$
\begin{aligned}
b_{hF}(\bm x) = L_f^{m}b(\bm x) +  S(b(\bm x)) + \alpha_m(\psi_{m-1}(\bm x)) \\+ L_gL_f^{m-1}b(\bm x)\bm u_{\min}\geq 0,
\end{aligned}
$$
and define $b_{hF}(\bm x)$ to be a HOCBF as in Def. \ref{def:hocbf}.

It is important to note that the relative degree of $b_{hF}(\bm x)$ with respect to dynamics (\ref{eqn:affine}) is only one, as we have $\psi_{m-1}(\bm x)$ in it. Thus,
we can get a feedback controller $K_{hF}(\bm x)$ that guarantees free conflict between the HOCBF constraint (\ref{eqn:hocbf}) and the control bounds (\ref{eqn:control}):
\begin{equation}\label{eqn:cbf_feaH}
\begin{aligned}
K_{hF}(\bm x) = \{\bm u\in \mathbb{R}^q: L_fb_{hF}(\bm x) + L_gb_{hF}(\bm x)\bm u \\+ \alpha_f(b_{hF}(\bm x))\geq 0\},
\end{aligned}
\end{equation}
if $b_{hF}(\bm x(0))\geq 0$, where $\alpha_f(\cdot)$ is a class $\mathcal{K}$ function.

\begin{theorem}\label{thm:feaH}
	If Problem \ref{prob:general} is initially feasible and the CBF constraint in (\ref{eqn:cbf_feaH}) corresponding to (\ref{eqn:feaH}) does not conflict with control bounds (\ref{eqn:control}) and (\ref{eqn:hocbf}) at the same time, any controller $\bm u \in K_{hf}(\bm x)$ guarantees the feasibility of problem (\ref{eqn:prob_qp}), subject to (\ref{eqn:prob_cbf})-(\ref{eqn:prob_ctrl}).
\end{theorem}
\textbf{Proof:} The proof is the same as Thm. \ref{thm:fea}.

Similar to the motivation for the analysis of the relative degree one case, we also  reformulate the constraint in (\ref{eqn:cbf_feaH}) in the form:
\begin{equation} \label{eqn:cbf_fea_reH}
\begin{aligned}
L_f^{m}b(\bm x) + L_gL_f^{m-1}b(\bm x)\bm u + S(b(\bm x)) + \alpha_m(\psi_{m-1}(\bm x)) \\+ \varphi(\bm x, \bm u) \geq 0.
\end{aligned}
\end{equation}
for some appropriate $\varphi(\bm x, \bm u)$. An obvious choice is $\varphi(\bm x, \bm u) = L_fb_{hF}(\bm x) + L_gb_{hF}(\bm x)\bm u + \alpha_f(b_{hF}(\bm x)) - L_f^{m}b(\bm x) - L_gL_f^{m-1}b(\bm x)\bm u - S(b(\bm x)) - \alpha_m(\psi_{m-1}(\bm x))$, which is a candidate function and we wish to simplify it. We define a set $U_{s}(\bm x)$ similar to (\ref{eqn:set_single}).

Similar to the last subsection, we just consider the case that the relative degree of $\varphi(\bm x, \bm u)$ is one, i.e., we have $\varphi(\bm x)$ from now on. Then, we have the following theorem to guarantee the feasibility of the problem (\ref{eqn:prob_qp}), subject to (\ref{eqn:prob_cbf})-(\ref{eqn:prob_ctrl}):
\begin{theorem} \label{thm:feasible_uniqueH}
	If $\varphi(\bm x)$ is a candidate function, $\varphi(\bm x(0))\geq 0, L_f\varphi(\bm x) \geq 0$, $L_g\varphi(\bm x) = \gamma L_gL_f^{m-1}b(\bm x)$, for some $\gamma > 0, \forall \bm x\in X$ and $\bm 0\in U$, then any controller $\bm u(t) \in U_s(\bm x), \forall t\geq 0$ guarantees the feasibility of the problem (\ref{eqn:prob_qp}), subject to (\ref{eqn:prob_cbf})-(\ref{eqn:prob_ctrl}).
\end{theorem}
\textbf{Proof:} The proof is the same as Thm. \ref{thm:feasible_unique}.

The approach to find a valid $\varphi(\bm x)$ is the same as the last subsection. The conditions in Thm. \ref{thm:feasible_uniqueH} are {\bf sufficient conditions} for the feasibility of the problem (\ref{eqn:prob_qp}), subject to (\ref{eqn:prob_cbf})-(\ref{eqn:prob_ctrl}).  Under the conditions in Thm \ref{thm:feasible_uniqueH}, we can also claim that $\varphi(\bm x) \geq 0$ is a single {\bf feasibility constraint} that guarantees the feasibility of the problem (\ref{eqn:prob_qp}), subject to (\ref{eqn:prob_cbf})-(\ref{eqn:prob_ctrl}) in the case that the safety constraint (\ref{eqn:safetycons}) is with high relative degree. We can get a feasible problem from the original problem (\ref{eqn:prob_qp}), subject to (\ref{eqn:prob_cbf})-(\ref{eqn:prob_ctrl}) in the form:
\begin{equation}\label{eqn:prob_qp2}
\min_{\bm u(t), \delta(t)} \int_{0}^T ||\bm u(t)||^2 + p\delta^2(t) dt
\end{equation}
subject to 
the feasibility constraint: (\ref{eqn:r0}) if the relative degree of $\varphi(\bm x, \bm u)$ is 0; otherwise, subject to the CBF constraint in (\ref{eqn:set_single}). The cost (\ref{eqn:prob_qp2}) is also subject to the HOCBF constraint (\ref{eqn:constraint}), the control bound (\ref{eqn:control}), and the CLF constraint:
\begin{equation} \label{eqn:prob_clf2}
L_{f}V(\bm x)+L_{g}V(\bm x)\bm u + \epsilon V(\bm x)\leq \delta(t),
\end{equation}
where $\varphi(\bm x)$ satisfies the conditions in Thm. \ref{thm:feasible_uniqueH} for (\ref{eqn:set_single}), and (\ref{eqn:r0}) is assumed to be non-conflicting with the HOCBF constraint (\ref{eqn:constraint}) and  the control bound (\ref{eqn:control}) at the same time.

\begin{remark}
	When we have multiple safety constraints, we can employ similar ideas to find sufficient conditions to guarantee problem feasibility. However, we also need to make sure that these sufficient conditions do not conflict with each other.
\end{remark}

\textbf{Example revisited.} We consider the example discussed in the beginning of this section, and demonstrate how we can find a single feasibility constraint $\varphi(\bm x(t))\geq 0$ for the ACC problem. It is obvious that $L_gL_fb(\bm x(t)) = -\frac{1}{M}$ in (\ref{eqn:safety_ex2}) does not change sign. The transformed control bound as in (\ref{eqn:ctrl_transH}) for (\ref{eqn:scons}) is
\begin{equation} \label{eqn:accc}
-c_dg\leq \frac{1}{M}u(t)\leq c_ag.
\end{equation}

The rewritten HOCBF constraint (\ref{eqn:safety_ex22}) can only conflict with the left inequality of (\ref{eqn:accc}). Thus, following (\ref{eqn:feaH}) and combining (\ref{eqn:safety_ex22}) with (\ref{eqn:accc}), the feasibility constraint is $b_{hF}(\bm x(t))\geq 0$, where
\begin{equation}\label{eqn:acc_fea}
\begin{aligned}
b_{hF}(\bm x(t)) =  \frac{F_r(v(t))}{M} + 2(p_1+p_2)(v_p - v(t)) \\+p_1p_2(z(t) - l_0) +c_dg.
\end{aligned}
\end{equation}

Since $\frac{F_r(v(t))}{M}\geq 0,\forall t\geq 0$, we can replace the last equation by
\begin{equation}\label{eqn:simplify}
\begin{aligned}
\hat b_{hF}(\bm x(t)) =  2(p_1+p_2)(v_p - v(t)) \\+p_1p_2(z(t) - l_0) +c_dg.
\end{aligned}
\end{equation}
The satisfaction of $\hat b_{hF}(\bm x(t))\geq 0$ implies the satisfaction of $ b_{hF}(\bm x(t))\geq 0$. Although the relative degree of (\ref{eqn:safety}) is two, the relative degree of $\hat b_{hF}(\bm x(t))$ is only one. We then define $\hat b_{hF}(\bm x(t))$ to be a CBF by choosing $\alpha_1(b(\bm x(t))) = kb(\bm x(t)), k>0$ in Def. \ref{def:hocbf}. Any control $u(t)$ should satisfy the CBF constraint (\ref{eqn:constraint}) which in this case is
\begin{equation} \label{eqn:fea_cbf}
\begin{aligned}
\frac{u(t)}{M} \leq \frac{F_r(v(t))}{M}  + (\frac{p_1p_2}{p_1+p_2}+k)(v_p - v(t)) \\+ \frac{kp_1p_2}{p_1+p_2}(z(t) - l_0) + \frac{kc_dg}{p_1 + p_2}
\end{aligned}
\end{equation}

In order to reformulate the last equation in the form of (\ref{eqn:cbf_fea_reH}), we try to find $k$ in the last equation. We require $\varphi(\bm x(t))$ to satisfy $L_g\varphi(\bm x(t)) \geq 0$ as shown in one of the conditions in Thm. \ref{thm:feasible_uniqueH}, thus, we wish to exclude the term $z(t)-l_0$ in $\varphi(\bm x(t))$ since its derivative $v_p - v(t)$ is usually negative. By equating the coefficients of the term $z(t) - l_0$ in (\ref{eqn:fea_cbf}) and (\ref{eqn:safety_ex22}), we have
\begin{equation}
\frac{kp_1p_2}{p_1+p_2} = p_1p_2
\end{equation}
Thus, we get $k = p_1 + p_2$. By substituting $k$ back into (\ref{eqn:fea_cbf}), we have
\begin{equation} \label{eqn:fea_cbf_re}
\begin{aligned}
\frac{u(t)}{M} \leq \frac{F_r(v(t))}{M}  + (p_1+p_2)(v_p - v(t)) \\+ {p_1p_2}(z(t) - l_0) + \varphi(\bm x(t))
\end{aligned}
\end{equation}
where 
\begin{equation} \label{eqn:acc_instance}
\varphi(\bm x(t)) = \frac{p_1p_2}{p_1+p_2}(v_p - v(t)) + c_dg
\end{equation}
It is easy to check that the relative degree of the last function is one, $L_f\varphi(\bm x(t)) = \frac{p_1p_2}{p_1+p_2}\frac{F_r(v(t))}{M} \geq 0$ and $L_g\varphi(\bm x(t)) = \frac{p_1p_2}{p_1+p_2}L_gL_fb(\bm x(t))$. Thus, all the conditions in Thm. \ref{thm:feasible_uniqueH} are satisfied except $\varphi(\bm x(0))\geq 0$ which depends on the initial state $\bm x(0)$ of system (\ref{eqn:veh}).  The single feasibility constraint $\varphi(\bm x(t))\geq 0$ for the ACC problem is actually a speed constraint (following (\ref{eqn:acc_instance})) in this case:
\begin{equation} \label{eqn:necess}
v(t)\leq v_p + \frac{c_dg(p_1 +p_2)}{p_1p_2}
\end{equation}
If $p_1 = p_2 = 1$ in (\ref{eqn:safety_ex22}), we require that the half speed difference between the front and ego vehicles should be greater than 
$-c_dg$ in order to guarantee the ACC problem feasibility.

We can find other sufficient conditions such that the ACC problem is guaranteed to be feasible by choosing different HOCBF definitions (different class $\mathcal{K}$ functions) in the above process.

\section{CASE STUDIES AND SIMULATIONS}

\label{sec:case}

In this section, we complete the ACC case study. All the computations and simulations were conducted in MATLAB. We used quadprog to solve the quadratic programs and ode45 to integrate the dynamics.

In addition to the dynamics (\ref{eqn:veh}), the safety constraint (\ref{eqn:safety}), the control bound (\ref{eqn:scons}), and the minimization of the cost $\int_{0}^{T}\left(\frac{u(t)-F_{r}(v(t))}{M}\right)^{2}dt$ introduced in Sec. \ref{sec:hocbf:moti}, we also consider a desired speed requirement $v\rightarrow v_d, v_d > 0$ in the ACC problem. We use the relaxed CLF as in (\ref{eqn:prob_clf}) to implement the desired speed requirement, i.e., we define a CLF $V = (v - v_d)^2$, and choose $c_1 = c_2 = 1, c_3 = \epsilon > 0$ in Def. \ref{def:clf}. Any control input should satisfy the CLF constraint (\ref{eqn:prob_clf}).

We consider the HOCBF constraint (\ref{eqn:safety_ex22}) to implement the safety constraint (\ref{eqn:safety}), and consider the sufficient condition (\ref{eqn:necess}) introduced in the last section to guarantee the feasibility of the ACC problem. We use a HOCBF with $m = 1$ to impose this condition, as introduced in (\ref{eqn:cbf_feaH}). We define $\alpha(\cdot)$ as a linear function in  (\ref{eqn:cbf_feaH}). 

Finally, we use the discretization method introduced in the end of Sec. \ref{sec:pre} to solve the ACC problem, i.e., 
We partition the time interval $[0,T]$ into a set of equal time intervals
$\{[0,\Delta t),[\Delta t,2\Delta t),\dots\}$, where $\Delta t>0$. In each
interval $[\omega\Delta t,(\omega+1)\Delta t)$ ($\omega=0,1,2,\dots$), we
assume the control is constant (i.e., the overall control will be piece-wise
constant), and reformulate the ACC problem as a sequence of QPs.
Specifically, at $t=\omega\Delta t$ ($\omega=0,1,2,\dots$), we solve 

\begin{equation} \label{eqn:objACC}
\bm u^*(t) = \arg\min_{\bm u(t)} \frac{1}{2}\bm u(t)^TH\bm u(t) + F^T\bm u(t)
\end{equation}
\[\begin{small}
\bm u(t)\! =\! \left[\begin{array}{c}  
\! u(t)\!\\
\!\delta(t)\!
\end{array} \right]\!,
H\! =\! \left[\begin{array}{cc} 
\frac{2}{M^2} & 0\\
0 & 2p_{acc}
\end{array} \right]\!, F\! =\!  \left[\begin{array}{c} 
\!\frac{-2F_r(v(t))}{M^2}\!\\
0
\end{array} \right]. \end{small}
\]
subject to
\[
A_{\text{clf}} \bm u(t) \leq b_{\text{clf}},
\]
\[
A_{\text{limit}} \bm u(t) \leq b_{\text{limit}},
\]
\[
A_{\text{hocbf\_safety}} \bm u(t) \leq b_{\text{hocbf\_safety}}, 
\]
\[
A_{\text{fea}} \bm u(t) \leq b_{\text{fea}}, 
\]
where $p_{acc} > 0$ and the constraint parameters are
$$
\begin{aligned}
A_{\text{clf}} &= [L_gV(\bm x(t)),\qquad -1],\\
b_{\text{clf}} &= -L_fV(\bm x(t)) - \epsilon V(\bm x(t)).
\end{aligned}
$$
$$
\begin{aligned}
A_{\text{limit}} &= \left[\begin{array}{cc} 
1, & 0\\
1, & 0
\end{array} \right],\\
b_{\text{limit}} &= \left[\begin{array}{c}  
c_aMg
-c_dMg
\end{array} \right].
\end{aligned}
$$
$$
A_{\text{hocbf\_safety}} = \left[\begin{array}{cc} 
\frac{1}{M}, & 0
\end{array} \right],
$$
$$
b_{\text{hocbf\_safety}} = \frac{F_r(v(t))}{M}  + (p_1+p_2)(v_p - v(t)) +p_1p_2(z(t) - l_0)
$$
$$
A_{\text{fea}} = \left[\begin{array}{cc} 
\frac{p_1p_2}{M(p_1+p_2)}, & 0
\end{array} \right],
$$
$$
b_{\text{fea}} = \frac{p_1p_2F_r(v(t))}{M(p_1+p_2)}  + \frac{p_1p_2}{p_1+p_2}(v_p - v(t)) + c_dg
$$

After solving (\ref{eqn:objACC}), we update (\ref{eqn:veh}) with $u^*(t)$, $\forall t\in (t_0 +\omega \Delta t, t_0 +(\omega+1) \Delta t)$.

\begin{table}
	\caption{Simulation parameters for the ACC problem}\label{table:param}
	\begin{center}
		\begin{tabular}{|c||c||c|c||c||c|}
			\hline
			Parameter & Value & Units &Parameter & Value & Units\\
			\hline
			\hline
			$v(0)$ & 6& $m/s$&	$z(0)$ & 100& $m$\\
			\hline
			$v_{p}$ & 13.89& $m/s$ & $v_d$ & 24& $m/s$\\
			\hline
			$M$ & 1650& $kg$ &g & 9.81& $m/s^2$\\
			\hline
			$f_0$ & 0.1& $N$ &$f_1$ & 5& $Ns/m$\\
			\hline
			$f_2$ & 0.25& $Ns^2/m$ &$l_0$ & 10& $m$\\
			\hline
			$\Delta t$ & 0.1& $s$&	$\epsilon$ & 10& unitless\\
			\hline
			$c_a(t)$ & 0.4& unitless&$c_d(t)$ & 0.4& unitless\\
			\hline
			$p_{acc}$ & 1& unitless &&&\\ 
			\hline
		\end{tabular}
	\end{center}
	
\end{table}

\begin{figure}[thpb]
	\centering
	\includegraphics[scale=0.475]{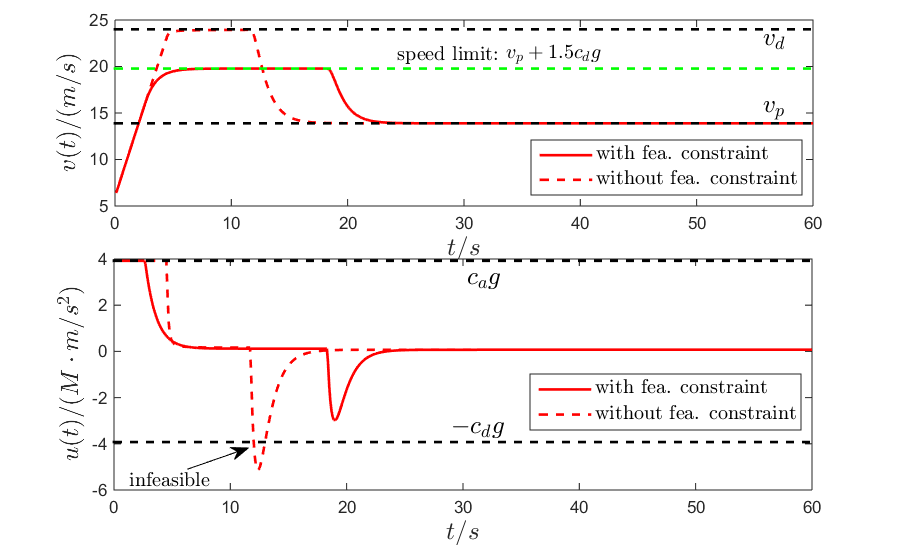}
	\caption{A simple case with $p_1 = 1, p_2 = 2$. The QP becomes infeasible when the ego vehicle exceeds the speed limit $v_p + 1.5c_dg$ from (\ref{eqn:necess}).}	
	\label{fig:single}
\end{figure}

\begin{figure}[thpb]
	\centering
	\includegraphics[scale=0.475]{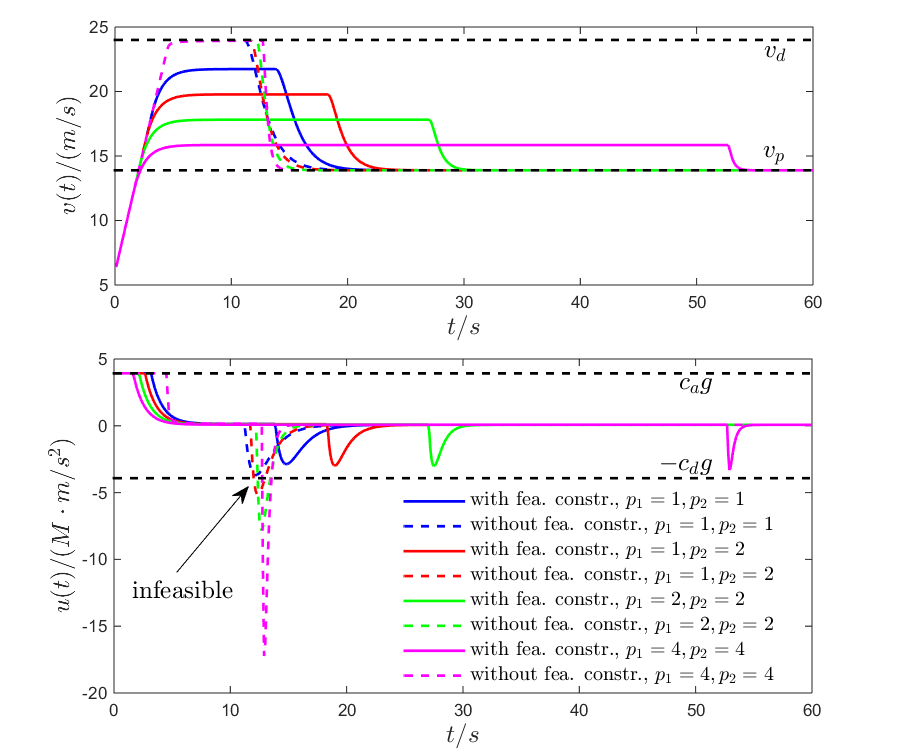}
	\caption{Speed and control profiles for the ego vehicle under different $p_1, p_2$, with and without feasibility condition (\ref{eqn:necess}).}	
	\label{fig:control_cmp}
\end{figure}

\begin{figure}[thpb]
	\centering
	\includegraphics[scale=0.5]{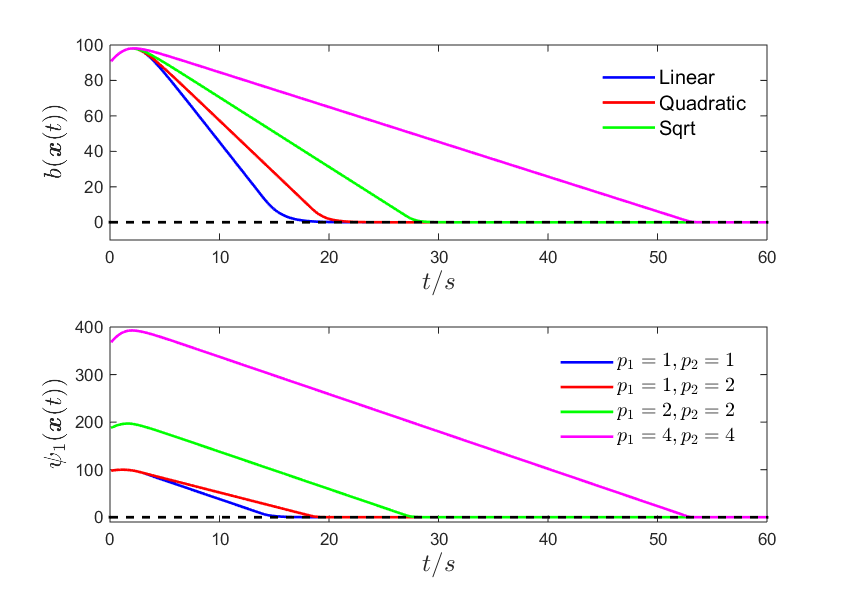}
	\caption{The variation of functions $b(\bm x(t))$ and $\psi_1(\bm x(t))$ under different $p_1, p_2$. $b(\bm x(t))\geq 0$ and $\psi_1(\bm x(t))\geq 0$ imply the forward invariance of the set $C_1\cap C_2$.}	
	\label{fig:set_cmp}
\end{figure}

The simulation parameters are listed in Table \ref{table:param}. We first present a case study in Fig. \ref{fig:single} showing that if the ego vehicle exceeds the speed constraint from the feasibility constraint (\ref{eqn:necess}), then the QP becomes infeasible. However, this infeasibility does not always hold since the feasibility constraint (\ref{eqn:necess}) is just a sufficient condition for the feasibility of QP (\ref{eqn:objACC}). In order to show how the feasibility constraint (\ref{eqn:necess}) can be adapted to different parameters $p_1, p_2$ in (\ref{eqn:safety_ex22}), we vary them and compare the solution without this feasibility sufficient condition in the simulation, as shown in Figs. \ref{fig:control_cmp} and \ref{fig:set_cmp}.

\begin{figure}[thpb]
	\centering
	\includegraphics[scale=0.52]{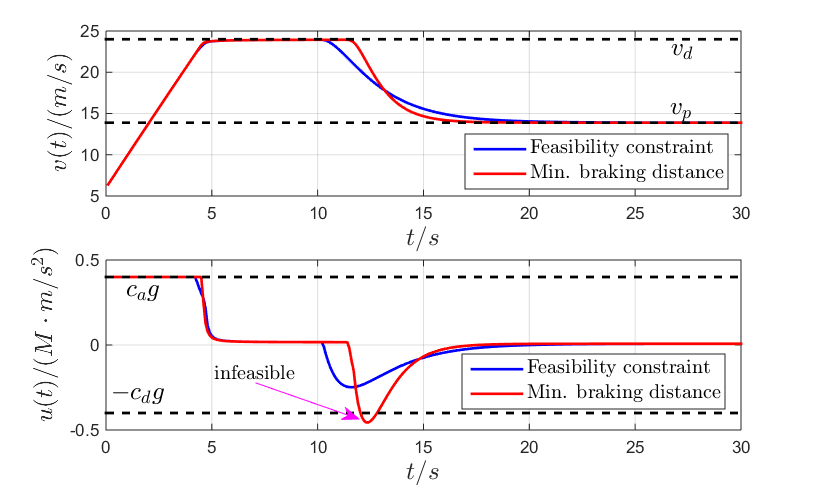}
	\caption{Comparison between the feasibility constraint (\ref{eqn:necess}) with $p_1 = 0.5, p_2 = 1$ and the minimum braking distance approach from \cite{Aaron2014}. The HOCBF constraint for (\ref{eqn:safety_bk}) in the minimum braking distance approach conflicts with the control bound (\ref{eqn:scons}).}	
	\label{fig:state}
\end{figure}

It follows from Figs. \ref{fig:control_cmp} and \ref{fig:set_cmp} that the QPs (\ref{eqn:objACC}) are always feasible with the feasibility constraint (\ref{eqn:necess}) under different $p_1, p_2$, while the QPs may become infeasible without this constraint. This validates the effectiveness of the feasibility constraint. We also notice that the ego vehicle cannot reach the desired speed $v_d$ with the feasibility condition (\ref{eqn:necess}); this is due to the fact that we are limiting the vehicle speed with  (\ref{eqn:necess}). In order to make the ego vehicle reach the desired speed, we choose $p_1, p_2$ such that the following constraint is satisfied.
\begin{equation}
v_p + c_dg \frac{(p_1 +p_2)}{p_1p_2}\geq v_d
\end{equation}
For example, the above constraint is satisfied when we select $p_1 = 0.5, p_2 = 1$ in this case. Then, the ego can reach the desired speed $v_d$, as the blue curves shown in Fig. \ref{fig:state}.

We also compare the feasibility constraint (\ref{eqn:necess}) with the minimum braking distance approach from \cite{Aaron2014}. This approach adds the minimum braking distance $\frac{0.5(v_p - v(t))^2}{c_dg}$ of the ego vehicle to the safety constraint (\ref{eqn:safety}):
\begin{equation}
\label{eqn:safety_bk}z(t) \geq\frac{0.5(v_p - v(t))^2}{c_dg} + l_{0},\forall t\geq0.
\end{equation}
Then, we can use a HOCBF with $m = 1$ (define $\alpha_1(\cdot)$ to be a linear function with slope 2 in Def. \ref{def:hocbf}) to enforce the above constraint whose relative degree is one. As shown in Fig. \ref{fig:state}, the HOCBF constraint for (\ref{eqn:safety_bk}) conflicts with the control bounds, and thus, the QP can still become infeasible.

\section{CONCLUSION \& FUTURE WORK}

\label{sec:conclusion}

We provide provably correct sufficient conditions for feasibility guarantee of constrained optimal control problems in this paper. These conditions are found by the proposed feasibility constraint method. We have demonstrated the effectiveness of sufficient feasibility conditions by applying them
to an adaptive cruise control problem. In the future, we will study the derivation of the necessary conditions of feasibility guarantee for constrained optimal control problems, or find less conservative sufficient conditions for specific dynamics. We will also try to figure out how to quickly find a single feasibility constraint for specific dynamics.

\bibliographystyle{IEEEtran}
\bibliography{CBF}

\begin{thebibliography}{10}
\providecommand{\url}[1]{#1}
\csname url@samestyle\endcsname
\providecommand{\newblock}{\relax}
\providecommand{\bibinfo}[2]{#2}
\providecommand{\BIBentrySTDinterwordspacing}{\spaceskip=0pt\relax}
\providecommand{\BIBentryALTinterwordstretchfactor}{4}
\providecommand{\BIBentryALTinterwordspacing}{\spaceskip=\fontdimen2\font plus
\BIBentryALTinterwordstretchfactor\fontdimen3\font minus
  \fontdimen4\font\relax}
\providecommand{\BIBforeignlanguage}[2]{{%
\expandafter\ifx\csname l@#1\endcsname\relax
\typeout{** WARNING: IEEEtran.bst: No hyphenation pattern has been}%
\typeout{** loaded for the language `#1'. Using the pattern for}%
\typeout{** the default language instead.}%
\else
\language=\csname l@#1\endcsname
\fi
#2}}
\providecommand{\BIBdecl}{\relax}
\BIBdecl

\bibitem{Bryson1969}
Bryson and Ho, \emph{Applied Optimal Control}.\hskip 1em plus 0.5em minus
  0.4em\relax Waltham, MA: Ginn Blaisdell, 1969.

\bibitem{Denardo2003}
E.~V. Denardo, \emph{Dynamic Programming: Models and Applications}.\hskip 1em
  plus 0.5em minus 0.4em\relax Dover Publications, 2003.

\bibitem{Rawlings2018}
J.~B. Rawlings, D.~Q. Mayne, and M.~M. Diehl, \emph{Model Predictive Control:
  Theory, Computation, and Design}.\hskip 1em plus 0.5em minus 0.4em\relax Nob
  Hill Publishing.

\bibitem{Aaron2014}
A.~D. Ames, J.~W. Grizzle, and P.~Tabuada, ``Control barrier function based
  quadratic programs with application to adaptive cruise control,'' in
  \emph{Proc. of 53rd IEEE Conference on Decision and Control}, 2014, pp.
  6271--6278.

\bibitem{Glotfelter2017}
P.~Glotfelter, J.~Cortes, and M.~Egerstedt, ``Nonsmooth barrier functions with
  applications to multi-robot systems,'' \emph{IEEE control systems letters},
  vol.~1, no.~2, pp. 310--315, 2017.

\bibitem{Xiao2019}
W.~Xiao and C.~Belta, ``Control barrier functions for systems with high
  relative degree,'' in \emph{Proc. of 58th IEEE Conference on Decision and
  Control}, Nice, France, 2019, pp. 474--479.

\bibitem{Tee2009}
K.~P. Tee, S.~S. Ge, and E.~H. Tay, ``Barrier lyapunov functions for the
  control of output-constrained nonlinear systems,'' \emph{Automatica},
  vol.~45, no.~4, pp. 918--927, 2009.

\bibitem{Wieland2007}
P.~Wieland and F.~Allgower, ``Constructive safety using control barrier
  functions,'' in \emph{Proc. of 7th IFAC Symposium on Nonlinear Control
  System}, 2007.

\bibitem{Boyd2004}
S.~P. Boyd and L.~Vandenberghe, \emph{Convex optimization}.\hskip 1em plus
  0.5em minus 0.4em\relax New York: Cambridge university press, 2004.

\bibitem{Aubin2009}
J.~P. Aubin, \emph{Viability theory}.\hskip 1em plus 0.5em minus 0.4em\relax
  Springer, 2009.

\bibitem{Prajna2007}
S.~Prajna, A.~Jadbabaie, and G.~J. Pappas, ``A framework for worst-case and
  stochastic safety verification using barrier certificates,'' \emph{IEEE
  Transactions on Automatic Control}, vol.~52, no.~8, pp. 1415--1428, 2007.

\bibitem{Wisniewski2013}
R.~Wisniewski and C.~Sloth, ``Converse barrier certificate theorem,'' in
  \emph{Proc. of 52nd IEEE Conference on Decision and Control}, Florence,
  Italy, 2013, pp. 4713--4718.

\bibitem{Panagou2013}
D.~Panagou, D.~M. Stipanovic, and P.~G. Voulgaris, ``Multi-objective control
  for multi-agent systems using lyapunov-like barrier functions,'' in
  \emph{Proc. of 52nd IEEE Conference on Decision and Control}, Florence,
  Italy, 2013, pp. 1478--1483.

\bibitem{Lindemann2018}
L.~Lindemann and D.~V. Dimarogonas, ``Control barrier functions for signal
  temporal logic tasks,'' \emph{IEEE Control Systems Letters}, vol.~3, no.~1,
  pp. 96--101, 2019.

\bibitem{Hsu2015}
S.~C. Hsu, X.~Xu, and A.~D. Ames, ``Control barrier function based quadratic
  programs with application to bipedal robotic walking,'' in \emph{Proc. of the
  American Control Conference}, 2015, pp. 4542--4548.

\bibitem{Wu2015}
G.~Wu and K.~Sreenath, ``Safety-critical and constrained geometric control
  synthesis using control lyapunov and control barrier functions for systems
  evolving on manifolds,'' in \emph{Proc. of the American Control Conference},
  2015, pp. 2038--2044.

\bibitem{Nguyen2016}
Q.~Nguyen and K.~Sreenath, ``Exponential control barrier functions for
  enforcing high relative-degree safety-critical constraints,'' in \emph{Proc.
  of the American Control Conference}, 2016, pp. 322--328.

\bibitem{Aaron2012}
A.~D. Ames, K.~Galloway, and J.~W. Grizzle, ``Control lyapunov functions and
  hybrid zero dynamics,'' in \emph{Proc. of 51rd IEEE Conference on Decision
  and Control}, 2012, pp. 6837--6842.

\bibitem{Galloway2013}
K.~Galloway, K.~Sreenath, A.~D. Ames, and J.~Grizzle, ``Torque saturation in
  bipedal robotic walking through control lyapunov function based quadratic
  programs,'' \emph{preprint arXiv:1302.7314}, 2013.

\bibitem{Xiao2020}
W.~Xiao, C.~Belta, and C.~G. Cassandras, ``Adaptive control barrier functions
  for safety-critical systems,'' in \emph{preprint in arXiv:2002.04577}, 2020.

\bibitem{Khalil2002}
H.~K. Khalil, \emph{Nonlinear Systems}.\hskip 1em plus 0.5em minus 0.4em\relax
  Prentice Hall, third edition, 2002.

\bibitem{Yang2019}
G.~Yang, C.~Belta, and R.~Tron, ``Self-triggered control for safety critical
  systems using control barrier functions,'' in \emph{Proc. of the American
  Control Conference}, 2019, pp. 4454--4459.

\bibitem{Xiao2020CDC}
W.~Xiao, C.~Belta, and C.~G. Cassandras, ``Feasibility guided learning for
  robust control in constrained optimal control problems,'' in \emph{to appear
  in CDC20, preprint in arXiv:1912.04066}, 2019.

\end{thebibliography}

\end{document}